\documentstyle{article}

\newcommand{\R}{{\bf R}}

\newcommand{\Z}{{\bf Z}}

\newcommand{\CP}{{\bf CP}}

\newcommand{\id}{{\rm id}}

\newcommand{\Ham}{{\rm Ham}}

\newcommand{\MS}{{\medskip}}

\newcommand{\NI}{{\noindent}}

\newcommand{\QED}{\hfill$\Box$\medskip}

\begin{document}

\title{Hofer's diameter and Lagrangian intersections}
\author{
Leonid
Polterovich
\\
Tel-Aviv University \\ (polterov@math.tau.ac.il)}

\date{January 12, 1998} 
\maketitle

\NI

The purpose of the present note is to show that the
group of Hamiltonian diffeomorphisms of the 2-sphere
has infinite diameter with respect to Hofer's metric.
For surfaces of higher genus this fact follows from
the energy-capacity inequality in the universal cover
(see Lalonde-McDuff [LM], and a recent preprint by M. Schwarz [S] 
which establishes
infiniteness of Hofer's diameter for arbitrary aspherical symplectic
manifolds).
Our approach is based on a reduction of the statement
to certain Lagrangian intersection problem, which in
the case of $S^2$ can be easily solved by existing methods.
The same approach works in some higher dimensional
situations (see discussion in section 6 below).

\medskip

{\bf 1.} Let $(M,\Omega)$ be a closed symplectic manifold.
Denote by $\cal F$ the space of all smooth functions
$F:M \times S^1 \to \R$ which satisfy a normalization condition
$\int_M F(x,t) \Omega^n = 0$ for every $t \in S^1$.
We write $\Ham (M,\Omega)$ for the group of Hamiltonian
diffeomorphisms of $M$. Given a transformation $\phi \in \Ham(M,\Omega)$,
its Hofer's norm $\rho(\id,\phi)$ is defined as
$$\inf \int_0^1 \max_x |F(x,t)| dt,$$
where the infimum is taken over all Hamiltonians $F \in \cal F$ which
generate $\phi$ (see [H]
\footnote{In [H] and [LM] the Hofer's norm is defined as
$\inf \int_0^1 \max_x F(x,t) - \min_x F(x,t) dt$. This norm
is equivalent to the one we consider, and the results 1.A, 1.B and 6.A
below remain valid for it as well.}
).

Consider the 2-sphere $S^2$ endowed with a symplectic form, and denote
by $L \subset S^2$ an equator (that is a simple closed curve which
divides the sphere into two parts of equal areas).

\proclaim Theorem 1.A. Let $\phi$ be a Hamiltonian diffeomorphism of 
$S^2$ generated by a Hamiltonian function $F \in \cal F$. Assume
that for some positive number $c$ holds $F(x,t) \geq c$ for
all $x \in L$ and $t \in S^1$. Then $\rho(\id,\phi) \geq c$.

\medskip

\proclaim Corollary 1.B. Hofer's diameter of $\Ham(S^2)$ is infinite.

\medskip

\NI
Indeed, choose an autonomous normalized Hamiltonian $F$ on $S^2$
which equals on the equator to an arbitrarily large constant.
Theorem 1.A implies that the corresponding Hamiltonian diffeomorphism
lies on an arbitrarily large distance from the identity.

The proof of 1.A is divided into several steps.

\medskip

{\bf 2.} Define a norm $||F|| =\max_{x,t} |F(x,t)|$ on $\cal F$.
First we note that the Hofer's distance can be expressed via this norm
as follows.

\proclaim Lemma 2.A. For every $\phi \in \Ham(M)$
$$\rho(\id,\phi) = \inf ||F||,$$ where the infimum is taken
over all Hamiltonians from $\cal F$ which generate $\phi$. 

\medskip

\NI
In the terminology of [P2]
this means that the "coarse" Hofer's norm
coincides with the usual one. The proof is based on a suitable choice
of the time reparametrizations. We present this elementary argument at the
end of the note.

\medskip

{\bf 3.} Denote by ${\cal H} \subset \cal F$ the subset of all
Hamiltonians from $\cal F$ which generate the identity map
(or in other words which correspond to loops of Hamiltonian diffeomorphisms).

\proclaim Lemma 3.A. Assume that a Hamiltonian diffeomorphism $\phi$
is generated by some $F \in \cal F$. Then 
$$\rho(\id,\phi) = \inf_{H \in {\cal H}} ||F-H||. $$

\medskip

\NI
{\it Proof:} Let $g =\{g_t\}, t \in [0;1]$ be any other path of Hamiltonian
diffeomorphisms joining the identity with $\phi$. Then $g_t = h_t \circ f_t$
for some loop $h =\{h_t\}$. Write $G$ and $H$ for the normalized Hamiltonians
of the paths $g$ and $h$ respectively. Then
$$G(x,t) = H(x,t) + F(h_t^{-1}x,t).$$
Assume that $G \in \cal F$. This is equivalent to the fact that $H \in \cal H$.
Set $H'(x,t) = -H(h_tx,t)$. Note that $H'$ generates a loop $\{h_t^{-1}\}$,
and thus belongs to $\cal H$. On the other hand the expression for $G$ above
implies that $||G|| = ||F - H'||$. The required statement follows immediately
from Lemma 2.A.
\QED

\medskip

Assume now that $M = S^2$, and $L$ is an equator.

\proclaim Lemma 3.B. For every $H \in \cal H$ there exists a point
$(x,t) \in L \times S^1$ such that $H(x,t) = 0$.

\medskip

\NI
{\it Proof of Theorem 1.A:} In view of 3.B the inequality $||F-H|| \geq c$
holds for all $H \in \cal H$. Thus 3.A implies that $\rho(\id,\phi) \geq c$.
This completes the proof.
\QED

\medskip

{\bf 4.} {\it Proof of 3.B:}

In order to prove Lemma 3.B recall the Lagrangian suspension
construction (which was used for study of Hofer's geometry in [P1]).
Let $L \subset M$ be a Lagrangian submanifold, and
$h = \{h_t\}$ be a loop of Hamiltonian diffeomorphisms. Then the
embedding $$L \times S^1 \to (M \times T^*S^1, \Omega \oplus dr \wedge dt)$$
which takes $(x,t)$ to $(h_tx, -H(h_tx,t),t)$ is Lagrangian. Denote its
image by $N(L,h)$. Note that homotopic loops $h$ lead to Hamiltonian
isotopic Lagrangian submanifolds $N(L,h)$ (this follows from a result
due to Weinstein [W]). 

Return now to the case when $M = S^2$ and $L$ is an equator.
We claim that for every loop $h$ the corresponding suspension
$N(L,h)$ is Hamiltonian isotopic to the "standard" Lagrangian
torus $N_0 = L \times \{r=0\}$. Indeed, it is known that
$\pi_1(\Ham(S^2)) = \Z_2$. If $h$ is homotopic to the constant
loop $\id$ then the claim follows from the fact that $N_0 = N(L,\id)$.
If $h$ is not contractible then it is homotopic to a path $f$ which
is given by the rotation around the axis orthogonal to the equator.
The corresponding Hamiltonian vanishes on $L$, and again we see
that $N_0 = N(L,f)$. This completes the proof of the claim.

Our next claim is that the Lagrangian submanifold $N_0$ has
the Lagrangian intersection property. Namely, if $N$ is an arbitrary
Lagrangian submanifold which is Hamiltonian isotopic to
$N_0$ then $N$ intersects $N_0$. 
This fact could be considered as {\it a stabilization} of the
(obvious) Lagrangian intersection property of the equator in $S^2$.
The proof is given in the next section.

Now we are ready to finish off the proof of the Lemma.
Assume that $N = N(L,h)$ where $h$ is a loop of Hamiltonian
diffeomorphisms of $S^2$. We see that $N(L,h)$ intersects $N_0$
and thus there exist $x \in L$ and $t \in S^1$ such that
$h_tx \in L$ and $H(h_tx,t) = 0$. This
completes the proof of Lemma 3.B.
\QED

\medskip

{\bf 5.} Here we prove the Lagrangian intersection property
which was used in the previous section. The proof is based
on the Floer homology theory for monotone Lagrangian submanifolds
developed by Y.-G. Oh (see [O2, section 8] and references therein).
Let $(M,\Omega)$
be a tame symplectic manifold, and $L \subset M$ be a 
closed monotone Lagrangian
submanifold whose minimal Maslov number is at least $2$.
Let $L'$ be the image of $L$ under a generic Hamiltonian isotopy
such that $L$ and $L'$ intersect transversally. Choose a generic
$\Omega$-compatible almost complex structure $J$ on $M$.
It is proved by Oh  that in this case the Floer homology
$HF(L,L')$ over $\Z_2$ is well defined and (up to an isomorphism)
does not depend on the generic choices of $L'$, the Hamiltonian isotopy
and $J$. If $HF(L,L') \neq 0$ then $L$ has the Lagrangian
intersection property, in other words $L \cap L'$ is non-empty.
In particular, this is the case when $L$ is the equator of $S^2$,
or more generally when $L= RP^n \subset \CP^n$ (see [O1]).

 Let $Z \subset T^*S^1$ be the zero section, and let
$Z' \subset T^*S^1$ be  
the graph
of the function $r = {\rm cos} 2\pi t$.
It is easy to see  that Floer homology $HF(L,L')$ is stable in the following
sense: 
$$ (5.A) \;\;\;\;\;\;HF(L\times Z, L' \times Z') = 
HF(L,L') \otimes (\Z_2 \oplus \Z_2).$$
Let us explain this formula.
First of all notice that the left hand side is well defined
since $L \times S^1$ is again monotone with the same minimal Maslov
number. Further, 
$$C(L\times Z,L' \times Z') = C(L,L') \otimes C(Z,Z'),$$
and
$$\partial  (x \otimes y) = \partial  (x) 
\otimes y + x \otimes \partial  (y).$$
Here $(C(.,.),\partial )$ stands for the Floer complex,
where $M \times T^*S^1$ is endowed with the split
almost complex structure, and the Hamiltonian isotopy between $L \times Z$
and $L' \times Z'$ splits as well.
The expression for the differentials above reflects the fact that
every discrete Floer's gradient trajectory is the product of
a discrete gradient trajectory in one factor and a constant trajectory
in the other factor. Clearly $C(Z,Z') = \Z_2 \oplus \Z_2$, and the
corresponding differential vanishes. Thus we get formula 5.A.

Applying 5.A to the case when $L$ is the equator of $S^2$ and taking
into account that $HF(L,L')$ is non-trivial we
obtain the required intersection property for $L \times S^1$.

\medskip

{\bf 6. A generalization and discussion.} 
Here we discuss some applications of our method to
higher-dimensional manifolds.
Let $(M,\Omega)$ be a closed symplectic
manifold which satisfies the following two conditions:

\MS

\NI
(i) The fundamental group of $\Ham(M,\Omega)$ is finite
(or, more generally, the Hofer's length spectrum of $\Ham (M,\Omega)$
is bounded away from infinity).

\medskip

\NI
(ii) There exists a monotone closed Lagrangian submanifold $L \subset M$
whose minimal Maslov number is at least 2, such that $HF(L,L') \neq 0$.

\MS

\proclaim Theorem 6.A. The Hofer's diameter of $\Ham(M,\Omega)$ is infinite.

\medskip

Note that this theorem applies to $S^2 \times S^2$ endowed with
the split symplectic form such that the areas of the factors are equal,
as well as to $\CP^2$. Indeed, take $L$ as the product of equators
in the first case, and as $RP^2$ in the second case, and note
that the corresponding fundamental groups are equal to $\Z_2 \oplus \Z_2$
and $\Z_3$ respectively. 

\MS

\NI
{\it Proof of 6.A:}
Let $F \in \cal F$ be a normalized autonomous Hamiltonian 
which equals to some positive constant $c$ on $L$.
Denote by $f = \{f_t\}, \; t \in [0;1]$ the corresponding
Hamiltonian isotopy which generates a Hamiltonian diffeomorphism
$\phi = f_1$.
Set $D=\Ham(M,\Omega)$ and let $({\tilde D}, \tilde \rho)$
be the universal cover of $D$ endowed with the lift
of the Hofer's metric. Here we take the identity element of $D$ as
the base point, and consider smooth paths and smooth homotopies
only.

Consider the lift $\psi$ of $\phi$
associated to the path $f$. Arguing exactly as in the proof of Theorem 1.A
we get that ${\tilde \rho}(\id, \psi) \geq c$.
Indeed, let $h$ be any contractible loop of Hamiltonian
diffeomorphisms. The condition (ii) combined with 5.A guarantees 
that the Lagrangian suspension $N(L,h)$ intersects $L \times Z$
in $M \times T^*S^1$. Then Lemmas 2.A and 3.A complete the job.

Denote by $u$ the upper
bound for the length spectrum of $G$. In other words, every closed loop
on $G$ can be homotoped to a loop whose Hofer's length does not exceed $u$.
The existence of such $u$ follows from the condition (i) above.
Since the Hofer's length structure is biinvariant, it follows that
for an arbitrary lift $\theta$ of $\phi$ holds
${\tilde \rho}(\psi, \theta) \leq u$. Thus the triangle inequality
for $\tilde \rho$ implies that ${\tilde \rho}(\id, \theta) \geq c-u$, and hence
$\rho(\id,\phi) \geq c-u$. Taking $c$ arbitrary large we get the
statement of the theorem.
\QED

\medskip

\NI
Note that Theorem 1.A is a particular case of Theorem 6.A, however
the proof of 1.A is in a sense more explicit than the one we just
completed . Namely it uses a more detailed
information about the Hamiltonians generating representatives	
of {\it all} homotopy classes of loops in $\Ham(S^2)$.
One can try to generalize this way of arguing to some
symplectic manifolds which violate condition (i) above.
For instance let $(M,\Omega)$ be the blow up of $\CP^2$
at one point endowed with the monotone symplectic structure.
It is proved recently in [AM] that the fundamental group
of $\Ham(M,\Omega)$ equals to $\Z$ and is generated by an
explicitly known $S^1$-action . Moreover the length spectrum is
unbounded (see [P2]). Denote by $H$ the normalized Hamiltonian of this
loop (see [P2] for the precise expression). It is easy to 
see that the zero level $\{H=0\}$ contains a Lagrangian torus $L$
which turns out to be {\it non-monotone}. It would be interesting
to check whether in spite of that non-monotonicity $L$ has the 
stable Lagrangian intersection property.
If it really does, then one concludes that the Hofer's diameter
in this case is infinite as well. 

In higher dimensions we face a difficulty. Namely, for 
applying our method we need an information about the fundamental group
of $\Ham(M,\Omega)$ which is not yet available. On the other hand
for some manifolds like $\CP^n$ or the product of $n$ spheres of equal areas 
one can use our approach in order to show that there exist arbitrarily
long Hamiltonian paths which cannot be shortened with fixed endpoints
at least in some homotopy classes.

Note finally that our arguments show that
if there exists a subset $L$ in $M$ with the open complement
such that every function $H \in \cal H$ vanishes at some point of 
$L \times S^1$ then the Hofer's diameter is infinite.
We took $L$ to be a Lagrangian submanifold in order to apply 
the existing Lagrangian
intersection techniques, however it would be quite natural
to work with more general (say, open) subsets. Here is a 
question which sounds as a first step towards such a generalization.
Take $L$ to be the complement of a tiny ball in $M$, and define
$N$ as the product of $L$ with the zero section in $T^*S^1$.
Is it true that $N$ always intersects its image under an
arbitrary Hamiltonian isotopy of $M \times T^*S^1$? Obviously
this is the case for split Hamiltonian isotopies. Can one use
the folding construction due to Lalonde and McDuff [LM] for
general Hamiltonian isotopies?

\medskip

{\bf 7.}{\it Proof of 2.A:} For $\phi \in \Ham(M,\Omega)$ 
set  $r(\id,\phi) = \inf ||F||$ where $F$ runs over 
all Hamiltonians $F \in \cal F$ which generate $\phi$.
Clearly $r(\id,\phi) \geq \rho(\id,\phi)$.
Our task is to prove the converse inequality. Fix a positive number $\epsilon$.
Choose a path $f = \{f_t\}, t \in [0;1]$ of Hamiltonian diffeomorphisms
generated by some $F \in \cal F$
such that $f_0 = \id, f_1 = \phi$ and such that
$\int_0^1 m(t) \leq \rho(\id,\phi) + \epsilon$ where $m(t) = \max_x |F(x,t)|$.
We can always assume that $m(t)$ is strictly positive
(geometrically this means that $f$ is a regular curve).
Denote by $\cal C$ the space of all $C^1$-smooth
orientation preserving diffeomorphisms
of $S^1$ which fix $0$. Note that for $a \in \cal C$
the path $ f_a =\{ f_{a(t)}\}$ is generated by the normalized Hamiltonian
function $F_a (x,t) = a'(t)F(x,a(t))$, where $a'$ denotes the derivative
with respect to $t$. Take now $a(t)$ as the inverse
of $$b(t) = {{\int_0^t m(s) ds}\over{\int_0^1 m(s)ds}}.$$
Note that 
$$||F_a|| = \max_t a'(t)m(a(t)) = \max_t (m(t)/b'(t)) = \int_0^1 m(t).$$
We conclude that $||F_a|| \leq \rho(\id, \phi) + \epsilon$.
Approximating $a$  in the $C^1$-topology by a smooth diffeomorphism from 
$\cal C$
we see that one can find a smooth normalized Hamiltonian, say $\tilde F$,
which generates $\phi$ and such that 
$||\tilde F|| \leq \rho(\id,\phi)+2\epsilon$. Since this can be done
for an arbitrary $\epsilon$ we conclude that $r(\id,\phi) \leq
\rho(\id,\phi)$. This completes the proof.
\QED

\bigskip

{\bf Acknowledgements.} This note was written during my stay
at ETH, Z\"urich. I thank the Forschungsinstitut f\"ur
Mathematik for the hospitality.

\bigskip

\NI
{\bf References}
\bigskip

\NI
[AM]
M.~Abreu and D.~McDuff, in preparation.
\MS

\NI
[H]
H.~Hofer, On the topological properties of symplectic maps,
{\it Proc. Roy. Soc. Edinburgh Sect A} {\bf 115}(1990), 25-38.
\MS

\NI
[LM]
F.~Lalonde and D.~McDuff, Hofer's
$L^{\infty}$-geometry: energy and stability of flows I,II,
{\it Invent. Math.} {\bf 122}  (1995), 1-69.
\MS

\NI
[O1] Y.-G.~Oh, Floer homology of Lagrangian intersections
and pseudo-holomor- phic discs I,II, {\it Comm. Pure Appl. Math.}
{\bf 46} (1993), 949-1012.
\MS

\NI
[O2] Y.-G.~Oh, Relative Floer and quantum cohomology and the
symplectic topology of Lagrangian submanifolds, in 
{\it Contact and Symplectic Geometry}, C.B. Thomas ed.,
Cambridge University Press 1996, pp. 201-267.
\MS

\NI
[P1] L.~Polterovich, Symplectic displacement energy
for Lagrangian submanifolds, {\it  Ergodic Th. and Dynam. Syst.}
{\bf 13} (1993), 357-367.
\MS

\NI
[P2]
L.~Polterovich, Hamiltonian loops and Arnold's 
principle, {\it Amer. Math. Soc. Transl.} (2) {\bf 180} (1997), 181-187.
\MS

\NI
[S] M.~Schwarz, A capacity for closed symplectically
aspherical manifolds, Prep- $\;\;\;$ rint 1997.
\MS

\NI
[W] A.~Weinstein, Cohomology of symplectomorphism groups and
critical values of Hamiltonians, {\it Math. Z.} {\bf 201} (1989), 75-82.

\end{document}